\documentclass [11pt, reqno]{amsart}
\usepackage{amsfonts,amssymb,amscd,amsmath,bbm}
\usepackage{enumerate,verbatim,calc,graphics,graphicx}
\usepackage{url}
\newtheorem{theorem}{Theorem}

\newtheorem*{rem}{Remark}
\begin{document}

\title [Evaluation of  Euler-phi sums of residue classes]
{\bfseries\large{Asymptotic evaluation of Euler-phi sums of various residue classes}}

\date{\today}

\author {Amrik Singh Nimbran}
\address{B3-304, Palm Grove Heights, Ardee City, Gurgaon, Haryana, INDIA }
\email{amrikn622@gmail.com}


\begin{abstract}
This note contains some asymptotic formulas for the sums of various residue classes of Euler's $\phi$-function.
\end{abstract}

\maketitle

\vspace{-0.7 cm}

\section{Introduction}

The \emph{phi-function} was introduced by Euler in connection with his generalization of Fermat's Theorem. It occurs without the functional notation in his 1759 paper \emph{Theoremata arithmetica nova methodo demonstrata} \cite{euler271}. In {\S 3} of his 1775 paper \cite{euler564}, Euler denotes by $\pi D$ ``the multitude of numbers less than D, and which have no common divisor with it" and then provides a table of $\pi D$ for $D=1$ to 100 writing $\pi 1=0.$ Gauss introduced the symbol $\phi$ in $\S 38$ of his \emph{Disquitiones Arithmeticae}(1801) with $\phi(1)=1.$ The function $\phi(n)$ denotes the number of positive integers not exceeding $n$ which are relatively prime to $n.$ Clearly, for $p$ prime, we have $\phi(p)=p-1.$

As Euler observed (Theorem 3, pp.81--82), if $p$ is a prime, the positive integers $\le p^k$ that are not relatively prime to $p^k$ are the $p^{k-1}$ multiples of $p: \; p, 2p, 3p, \ldots, p^{k-1}\cdot p.$ So $\phi(p^k)=p^{k}-p^{k-1}=p^{k}(1 -\frac{1}{p})=p^{k-1}(p-1),$ and $\sum_{j=0}^{k} \phi(p^j)=(p-1)[1+p+p^2+\dots+p^{k-1}]=p^k.$ Furthermore, if $(a, b)=1,$ then $\phi(a \, b)= \phi(a) \, \phi(b).$ Thus if $m$ has the prime factorization $m=p_1^{r_1} p_2^{r_2} \cdots p_k^{r_k},$ then $\phi(m) =p_1^{r_1-1} p_2^{r_2-1} \cdots p_k^{r_k-1} (p_1-1) (p_2-1)\cdots (p_k-1).$ And, $\phi(m^k)=m^{k-1} \phi(m).$ Also, if $(a, b)=d,$ then $\displaystyle \phi(a \, b)=\phi(a) \, \phi(b) \frac{d}{\phi(d)}.$ As Gauss showed: $\displaystyle \sum_{d|n} \phi(d)=\sum \phi(n/d)=n.$

The value of $\phi(n)$ fluctuates as $n$ varies. Since averages sooth out fluctuations, it may be fruitful to study the arithmetic mean $\displaystyle \frac{\Phi(n)}{n},$ where $\Phi(n)=\sum_{m=1}^n \phi(m).$

In 1874, Mertens obtained \cite[p.122]{dickson}\cite{mertens} an asymptotic value for $\Phi(N)$ for large $N.$ He employed the function $\mu(n)$ and proved that
\[
\sum_{m=1}^{G} \phi(m)=\frac{1}{2} \sum_{n=1}^{G} \mu (n) \left\lbrace \left[\frac{G}{n} \right]^2 +\left[\frac{G}{n} \right]\right\rbrace =\frac{3}{\pi^2} G^2+\Delta
\]
with $|\Delta |<G(\frac{1}{2} \ln G +\frac{1}{2} \gamma +\frac{5}{8}) +1,$ where $\gamma$ is Euler's constant and $\mu (n)$
is the M\"{o}bius function defined as
\[
\mu(n)=
\begin{cases}
1  & \text{if} \, n=1, \\
(-1)^r & \text{if $n$ is product of $r$ distinct prime numbers}, \\
0  & \text{if $n$ has one or more repeated prime factors}. \\
\end{cases}
\]

If $(a, \, b)=1, \; \mu(a \, b)=\mu(a) \, \mu(b).$ Further, $\displaystyle \sum_{d|n} \mu(d)=0 \quad (n>1).$

For any positive integer $n,$ we have\cite[pp.78--80]{andrews}:
\[
\phi(n)=\sum_{d|n} \frac{n}{d} \, \mu(d)=\sum_{d|n} d \, \mu\left(\frac{n}{d} \right).
\]

It is shown in \cite[p.268 Theorem 330]{hardy}\cite[pp.61-62]{apostol} that:
\begin{equation}
\Phi(n)=\frac{3 n^2}{\pi^2} +O(n \, \ln n). \label {eq:2}
\end{equation}

To prove (\ref{eq:2}), we  may recall here Euler's \emph{zeta function} and identity:
\[
\zeta(s)=\sum_{n=1}^\infty \frac{1}{n^s}=\prod_{p-prime} \left(1 -\frac{1}{p^s}\right)^{-1}, \; \Re(s)>1.
\]

Since for $s>1,$
\[
\frac{1}{\zeta(s)}=\prod_{p} \left(1 -p^{-s}\right)= \prod \lbrace 1 +\mu(p)p^{-s} +\mu(p^2)p^{-2s} +\ldots\rbrace=\sum_{n=1}^\infty \frac{\mu(n)} {n^{s}}
\]
and
\[
\phi(n)=n \sum_{d|n} \frac{\mu(d)}{d}
\]
we have:
\[
\Phi(n)=\sum_{m=1}^n \phi(m)=\sum_{m=1}^n m \sum_{d|m} \frac{\mu(d)}{d}=\sum_{d d^{\prime} \le n} d^{\prime} \mu(d)
=\sum_{d=1}^n \mu(d) \sum_{d^{\prime}=1}^{\lfloor\frac{n}{d}\rfloor} d^{\prime}.
\]
That is,
\[
\Phi(n)=\sum_{d=1}^n \mu(d) \left\lbrace \frac{1}{2} \left\lfloor\frac{n}{d}\right\rfloor \left( \left\lfloor\frac{n}{d}\right\rfloor +1\right) \right\rbrace =\frac{1}{2} \sum_{d=1}^n \mu(d) \left\lbrace \frac{n^2}{d^2} +O\left(\frac{n}{d}\right)\right\rbrace,
\]
leading to
\begin{align*}
\Phi(n)=&\frac{n^2}{2} \sum_{d=1}^n \frac{\mu(d)}{d^2}+O\left(n \sum_{d=1}^n \frac{1}{d}\right)\\
=& n^2\sum_{d=1}^\infty \frac{\mu(d)}{d^2}- n^2 \sum_{d=n+1}^\infty \frac{\mu(d)}{d^2} +O(n \ln n).\\
=&\frac{n^2}{2\zeta(2)} +O\left(n^2 \sum_{d=n+1}^\infty \frac{1}{d^2}\right) +O(n \ln n).
\end{align*}
Or,
\[
\Phi(n)=\frac{n^2}{2\zeta(2)} +O(n) +O(n \ln n)=\frac{3 n^2}{\pi^2} +O (n \ln n).
\]

Lehmer studied sums of $\phi(n)$ in \cite{lehmer1} and revisited in \cite{lehmer2}. I seek here an extension of Lehmer's formula occurring in \cite{lehmer2} by using his argument.

\section{Asymptotic summation of $\phi(pn)$}

Since $\phi(2^k)=2^{k-1},$ so: $\phi(4m+2)=\phi(2m+1); \quad \phi(4m)=2\phi(2m).$

Denoting $\displaystyle \Phi_{e}(n)=\sum_{m\le n; \, m \, even} \phi(m)$ and $\displaystyle\Phi_{o}(n)=\sum_{m\le n; \, m \, odd} \phi(m),$ and using the relation:
\[
\Phi_{e}(n)=\Phi_{o}(n/2) +2\Phi_{e}(n/2) =\Phi(n/2) +\Phi_{e}(n/2),
\]
Lehmer \cite{lehmer2} deduced: $\displaystyle \Phi_{e}(n)=\sum_{\lambda=1}^{\ell}  \Phi_{e}(n/2) \quad (\ell=[\ln n /\ln 2])$
and then used the formula (\ref{eq:2}) to derive:
\begin{equation}
\Phi_{e}(n)=\left(\frac{n}{\pi}\right)^2 +O(n \, \ln n); \quad \Phi_{o}(n)=2 \left(\frac{n}{\pi}\right)^2 +O(n \, \ln n). \label {eq:3}
\end{equation}

Let $\displaystyle \Phi_{\displaystyle r_i}(n)=\sum_{k=1}^m \phi(kp-i),$ with fixed $i=0, 1, 2, \dots, p-1$ and $(mp-i)\le n.$ Then
$\displaystyle \Phi_{r_0}(n)=(p-1) \sum_{i=1}^{p-1} \Phi_{r_i}(n/p) +p \, \Phi_{r_0} (n/p).$ Hence,

\[
\Phi_{\displaystyle r_0}(n)=(p-1) \, \Phi (n/p) +\Phi_{\displaystyle r_0} (n/p).
\]

Mimicking Lehmer's proof, we see that for any prime $p,$
\begin{align*}
\Phi_{\displaystyle r_0} (n) & = (p-1) \, 3 \pi^{-2} \, n^2 \, \sum_{\lambda=1}^q p^{-2\lambda} +O(n \log n) \quad (q = [\ln n / \ln p])\\
 & = \frac{3(p-1)}{p^2-1} \pi^{-2} \, n^2 + O \left(n^2 \int_q^\infty (p^{-2})^{t} \, {dt} \right) +O(n \log n) \\
 & = \frac{3}{p+1} \pi^{-2} \, n^2 + O(n \log n). \Box
\end{align*}

The last asymptotic formula implies the following theorem:
\begin{theorem}
For any prime $p,$ we have:
\begin{equation}
\lim_{m \to \infty} \frac{\displaystyle \sum_{k=1}^m \phi(pk)}{(pm)^2}=\frac{3}{(p+1)\pi^2}. \label {eq:4}
\end{equation}
\end{theorem}

If the set $\mathbb{N}$ is partitioned into $p$ residue classes modulo $p,$ we will have one class consisting of composite numbers of the form $pm$ while the remaining $p-1$ classes contain nearly an equal number of prime numbers, and the ratio of the cumulative sums of the two types of classes will be $p:(p-1).$ The rationale behind the first part of the statement is found in Dirichlet's famous theorem relating to primes in arithmetic progressions: \textit{every arithmetic progression, with the first member and the difference being coprime, will contain infinitely many primes.} In other words, if $k>1$ is an integer and $(k, \ell)=1,$ then there are infinitely many primes of the form $kn+\ell,$ where $n$ runs over the positive integers. If $k$ is a prime $p,$ then $\ell$ is one of the numbers $1, 2, \ldots, p-1.$

Let us recall here the arithmetic function known as the \emph{Mangoldt function} which is defined as:
\[
\Lambda(n)=
\begin{cases}
\ln p, \; \text{if} \; n=p^m \; \text{for some prime} \; p \; \text{and positive integer} \; m,\\
0 \; \text{otherwise}.
\end{cases}
\]
This function has an important role in elementary proofs of the prime number theorem which states that if $\pi(n)$ denotes the number of primes $\le n,$ then  $\displaystyle \pi (n) \sim \frac{n}{\ln n}.$ We have (\cite[pp.253-254]{hardy}) for $n\ge 1:$

\[
\Lambda(n)=\sum_{d|n} \mu \left(\frac{n}{d}\right) \, \ln d=\sum_{d|n} \mu(d) \, \ln \left(\frac{n}{d}\right)=-\sum_{d|n} \mu(d) \ln d,
\]
and
\[
\sum_{d|n} \Lambda(d)=\ln n.
\]

Further \cite[p.348]{hardy}\cite[p.89]{apostol},
\[
\sum_{n \le x} \frac{\Lambda(n)}{n}=\ln x + O(1),
\]
whence
\begin{equation}
\sum_{p \le x} \frac{\ln p}{p}=\ln x +O(1).  \label{eq:a}
\end{equation}

This related result is well-known\cite[p.148]{apostol}:
\begin{equation}
\sum_{\substack{p\le x;\\ p\equiv \ell \, (mod \, k)}} \frac{\ln p}{p}=\frac{1}{\phi(k)} \ln x +O(1), \label{eq:b}
\end{equation}
where the sum is extended over those primes $p \le x$ which are congruent to $\ell \; (\text{mod} \, k).$ Since $\ln x \to \infty$ as $x \to \infty$ this relation implies that there are infinitely many primes $p \equiv \ell (\text{mod} \, k),$ hence infinitely many in the progression $kn+\ell.$ Since the principal term on the right hand side in (\ref{eq:b}) is independent of $\ell,$ therefore it not only implies Dirichlet's theorem but it also shows \cite[p. 148]{apostol} that the primes in each of the $\phi (k)$ reduced residue classes $(\text{mod} \, k)$ make the same contribution to the principal term in (\ref{eq:a}), that is, the \textit{primes are equally distributed among $\phi (k)$ reduced residue classes $(\text{mod} \, k).$} We thus have a prime number theorem for arithmetic progressions \cite[p. 154]{apostol}: If $\pi_{\ell} (x)$ counts the number of primes $\le x$ in the progression $kn+\ell,$ then
\[
\pi_{\ell} (x) \sim \frac{\pi(x)}{\phi(k)} \sim \frac{1}{\phi(k)} \frac{x}{\ln x}.
\]

Hence as $m \to \infty, \; \Phi_{r_i}(m) \sim \Phi_{r_j}(m), \, i, j \ne 0.$ And so, we deduce from (\ref{eq:2}) and our Theorem 1 the following result:
\begin{theorem}
For any prime $p,$ we have for each $i=1, 2, 3, \dots, p-1,$
\begin{equation}
\lim_{m \to \infty} \frac{\displaystyle \sum_{k=1}^m \phi(pk-i)}{(pm)^2} =\frac{3p}{(p^{2} -1)\pi^2};  \;
\end{equation}
\end{theorem}

We will now obtain asymptotic evaluation of the sums of residue classes modulo $p$ for the $\phi$-function.

Since $\phi(4m-2)=\phi(2m-1); \; \phi(4m)=2\, \phi(2m)$ and as $n\to \infty, \; \Phi(2n-1)=\displaystyle \sum_{m=1}^n \phi(2m-1)=2 \, \Phi(2n)=2 \displaystyle \sum_{m=1}^n \phi(2m),$ so we have:

\[
\lim_{n \to \infty} \frac{\Phi(4n-2)}{(4n)^2} = \lim_{n \to \infty} \frac{\Phi(4n)}{(4n)^2} =\frac{1}{2\pi^2}.
\]

Further, as $n\to \infty; \; \Phi(2n-1)=\Phi(4n-3) +\Phi(4n-1)=2 \, \Phi(2n)=2 \, \Phi(4n-2) +2 \, \Phi(4n)$ and the two forms $4k-3, \; 4k-1$ yield almost equal number of primes, so we have:

\[
\lim_{n \to \infty} \frac{\Phi(4n-3)}{(4n)^2} = \lim_{n \to \infty} \frac{\Phi(4n-1)}{(4n)^2}=\frac{1}{\pi^2}.
\]

Again,
\[
\lim_{n \to \infty} \frac{\Phi(6n-4)}{(6n)^2} +\lim_{n \to \infty} \frac{\Phi(6n-2)}{(6n)^2} +\lim_{n \to \infty} \frac{\Phi(6n)}{(6n)^2}=\frac{1}{\pi^2}
\]
and
\[
\lim_{n \to \infty} \frac{\Phi(6n-5)}{(6n)^2} +\lim_{n \to \infty} \frac{\Phi(6n-3)}{(6n)^2} +\lim_{n \to \infty} \frac{\Phi(6n-1)}{(6n)^2}=\frac{2}{\pi^2}.
\]
Further
\[
\lim_{n \to \infty} \frac{\Phi(6n-4)}{(6n)^2} = \quad \lim_{n \to \infty} \frac{\Phi(6n-2)}{(6n)^2} =\frac{3}{2} \lim_{n \to \infty} \frac{\Phi(6n)}{(6n)^2}
\]
and
\[
\lim_{n \to \infty} \frac{\Phi(6n-5)}{(6n)^2} = \quad \lim_{n \to \infty} \frac{\Phi(6n-1)}{(6n)^2} =\frac{3}{2} \lim_{n \to \infty} \frac{\Phi(6n-3)}{(6n)^2}.
\]
Still further,
\[
\lim_{n \to \infty} \frac{\Phi(3(2n-1))}{(6n)^2} =2 \lim_{n \to \infty} \frac{\Phi(3(2n))}{(6n)^2}.
\]
So we deduce these results:
\[
\lim_{n \to \infty} \frac{\Phi(6n-4)}{(6n)^2} = \quad \lim_{n \to \infty} \frac{\Phi(6n-2)}{(6n)^2} =\frac{3}{8\pi^2}.
\]

\[
\lim_{n \to \infty} \frac{\Phi(6n-3)}{(6n)^2} =\frac{1}{2 \pi^2};  \quad \lim_{n \to \infty} \frac{\Phi(6n)}{(6n)^2} =\frac{1}{4\pi^2}.
\]

\[
\lim_{n \to \infty} \frac{\Phi(6n-5)}{(6n)^2} = \quad \lim_{n \to \infty} \frac{\Phi(6n-1)}{(6n)^2} =\frac{3}{4\pi^2}.
\]

In fact, we have following the general theorem based on two facts: (i) the sum of all odd residue classes equals twice the sum of all even classes, and (ii) the ratio of residue classes modulo $p$ containing primes to the class having only composite numbers is $\displaystyle \frac{p}{p-1}:1.$

\begin{theorem}
For an odd prime $p,$
\begin{align*}
\lim_{n \to \infty} \frac{\Phi(2pn -(2p-2))}{(2pn)^2}=\lim_{n \to \infty} \frac{\Phi(2pn -2)}{(2pn)^2}=\\
\lim_{n \to \infty} \frac{\Phi(2pn -(2p-4))}{(2pn)^2}=\lim_{n \to \infty} \frac{\Phi(2pn -4)}{(2pn)^2}=\\
\dots\\
\lim_{n \to \infty} \frac{\Phi(2pn -(p+1))}{(2pn)^2}=\lim_{n \to \infty} \frac{\Phi(2pn -(p-1))}{(2pn)^2}=
\frac{p}{(p^{2}-1)\pi^2};\\
\lim_{n \to \infty} \frac{\Phi(2pn)}{(2pn)^2}=\frac{1}{(p+1)\pi^2}.
 \end{align*}
And
\begin{align*}
\lim_{n \to \infty} \frac{\Phi(2pn -(2p-1))}{(2pn)^2}=\lim_{n \to \infty} \frac{\Phi(2pn -1)}{(2pn)^2}=\\
\lim_{n \to \infty} \frac{\Phi(2pn -(2p-3))}{(2pn)^2}=\lim_{n \to \infty} \frac{\Phi(2pn -3)}{(2pn)^2}=\\
\dots\\
\lim_{n \to \infty} \frac{\Phi(2pn -(p+2))}{(2pn)^2}=\lim_{n \to \infty} \frac{\Phi(2pn -(p-2))}{(2pn)^2}=\frac{2p}{(p^{2}-1)\pi^2};\\
\lim_{n \to \infty} \frac{\Phi(2pn -p)}{(2pn)^2}=\frac{2}{(p+1)\pi^2}.
 \end{align*}
\end{theorem}

\begin{rem}
If $m<p,$ then $m$ cannot divide $p.$ Also, $p$ cannot divide $2m$ and $2m-1$ simultaneously; it may not divide either. So gcd $(p, 2m)=1 \, \text{or} \; p$ and gcd $(p, 2m-1)=p \, \text{or} \, 1.$ Hence, $\phi(p \,(2m))=(p-1) \, \phi(2m) $ or $p \, \phi(2m);$ and $\phi(p (2m-1))=p \, \phi(2m-1)$ or $(p-1) \, \phi(2m-1)$ depending on $m.$ Lehmer proved that $\displaystyle \lim_{n \to \infty} \frac{\Phi(2n-1)}{\Phi(2n)}=2.$ Hence, $\displaystyle \lim_{n \to \infty} \frac{\sum_{m=1}^{n} \phi(p \,(2m-1))}{\sum_{m=1}^{n} \phi(p \,(2m))}=2.$
\end{rem}

\noindent {\bf Acknowledgement:} The author is thankful to Prof Paul Levrie for his helpful comments and the anonymous referee for his suggestions which made the presentation concise.

\end{document}